\documentclass[a4paper,fleqn]{cas-sc}
\usepackage{bm}
\usepackage{float}
\usepackage{placeins}
\usepackage[numbers]{natbib}
\usepackage{amsmath, amssymb, amsthm}
\usepackage{algorithm, algorithmic}
\usepackage{graphicx}
\usepackage{xcolor}
\usepackage{bbm}
\usepackage{subfigure}
\usepackage{multirow}
\usepackage{empheq}
\usepackage{kotex}
\usepackage{array}

\def\tsc#1{\csdef{#1}{\textsc{\lowercase{#1}}\xspace}}
\tsc{WGM}
\tsc{QE}
\tsc{EP}
\tsc{PMS}
\tsc{BEC}
\tsc{DE}

\newcolumntype{C}{>{\footnotesize}c}
\newcolumntype{R}{>{\footnotesize$}r<{$}}
\newcolumntype{L}{>{\footnotesize$}l<{$}}
\newcolumntype{P}{>{\footnotesize$}c<{$}}
\graphicspath{{./}}

\begin{document}
\let\WriteBookmarks\relax
\def\floatpagepagefraction{1}
\def\textpagefraction{.001}
\shorttitle{A Lower-envelope Eikonal Neural Operator}
\shortauthors{H. Son, S. W. Cho, and J. Hahn}

\title [mode = title]{A Lower-envelope Eikonal Neural Operator}       
\author[1]{Hwijae Son}
\ead{hwijaeson@konkuk.ac.kr}
\credit{Conceptualization, Data curation, Methodology, Validation, Software, Visualization, Writing -- original draft, Writing -- review and editing}
\affiliation[1]{organization={Department of Mathematics, Konkuk University},
                addressline={20 Neungdong-ro, Gwangjin-gu}, 
                city={Seoul},
                postcode={05029}, 
                country={Republic of Korea}}

\author[2]{Sung Woong Cho}
\ead{swcho95kr@inha.ac.kr}
\credit{Conceptualization, Methodology, Validation, Software, Writing -- original draft, Writing -- review and editing}
\affiliation[2]{organization={Department of Data Science, Inha University},
                addressline={313, Dokbae-ro, Michuhol-gu}, 
                postcode={22188}, 
                city={Incheon},
                country={Republic of Korea}}

\author[3]{Jooyoung Hahn}
\credit{Methodology, Supervision, Validation, Writing -- original draft, Writing -- review \& editing}
\affiliation[3]{organization={Department of Software Engineering, Faculty of Nuclear Sciences and Physical Engineering, Czech Technical University in Prague},
	addressline={Trojanova 339/13},
	city={Prague},
	postcode={120 00},
	country={Czech Republic}}
\cormark[1]
\ead{jooyoung.hahn@fjfi.cvut.cz}

\cortext[cor1]{Corresponding author}

\begin{abstract}
We propose the Lower-envelope Eikonal Neural Operator (LENO), a physics-informed neural operator for eikonal equations with heterogeneous velocity fields and general source sets. LENO encodes the variational minimum principle of viscosity solutions by constructing source-wise candidate fields and taking their pointwise minimum. This lower-envelope structure enforces the source condition exactly and represents nonsmooth multi-source solutions without a soft boundary penalty. The model is trained without precomputed solution data, using a Godunov upwind residual as the sole training signal. Experiments on point sources and curve sources, sample-dependent piecewise velocity fields, and OpenFWI CurveVel-A and Style-A media show accurate agreement with references of Fast Marching Method (FMM). Additional validation shows that the Godunov residual outperforms conventional finite difference residuals on a high-resolution test, that velocity conditioning and the lower-envelope transform are essential in realistic heterogeneous media, and that the trained operator supports efficient batched evaluation and resolution transfer.
\end{abstract}

\begin{graphicalabstract}
\includegraphics[width=\linewidth]{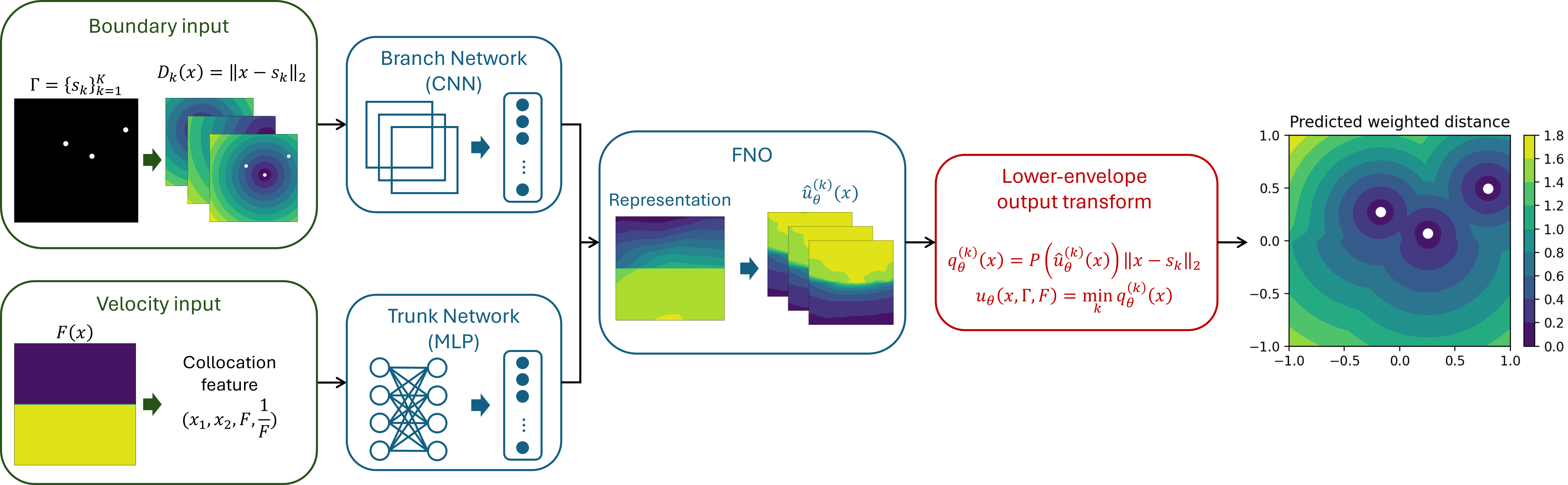}
\end{graphicalabstract}

\begin{highlights}
    \item We introduce LENO, a physics-informed neural operator that represents eikonal viscosity solutions as lower envelopes of source-wise candidate weighted distance fields.

    \item The lower-envelope transform enforces the source condition exactly and removes the need for a soft boundary penalty in residual-only training.

    \item A Godunov upwind residual provides a viscosity-consistent training signal and outperforms conventional finite difference residuals in a high-resolution stencil study.

    \item Experiments validate LENO on point sources, open and closed curve sources, piecewise velocity fields, and OpenFWI CurveVel-A and Style-A media.

    \item Ablation, runtime, and resolution-transfer studies confirm the roles of velocity conditioning, the lower-envelope transform, batched evaluation, and finer-grid inference.
\end{highlights}

\begin{keywords}
Eikonal equation \sep Viscosity solution \sep Physics-informed neural operator
\end{keywords}

\maketitle

\section{Introduction}\label{sec:intro}

We consider the eikonal equation on a bounded convex domain $\Omega\subset\mathbb{R}^2$,
\begin{empheq}[left=\empheqlbrace]{equation}
    \begin{aligned}\label{eq:eikonal_main}
		F(x)\,\|\nabla u(x)\| &= 1, && x\in\Omega\setminus\Gamma,\\
		u(x) &= 0, && x\in\Gamma,
	\end{aligned}
\end{empheq}
where $\|\cdot\|$ is the Euclidean norm, $F:\bar{\Omega}\to(0,\infty)$ is the velocity field (with reciprocal $1/F$, the slowness), and the source set $\Gamma$ is a nonempty compact subset of $\bar{\Omega}$. The eikonal equation is the prototypical static first-order Hamilton-Jacobi equation. The viscosity solution of~\eqref{eq:eikonal_main} carries two equivalent and concrete meanings at once. It measures the slowness-weighted distance from each point to the source set, and equally the first-arrival time at that point of a front that starts on the source set and spreads through the medium at the local velocity. This dual reading makes the same equation recur across disciplines. In seismology it produces the first-arrival traveltime fields used for tomography and earthquake source localization; the traveltime field is also the forward map inside adjoint-state tomography, where the velocity is reconstructed from observed arrivals~\cite{tong2021adjoint}. In image processing and computer vision it supplies distance transforms, geodesic distances, and the front used in segmentation and registration; on manifolds and triangulated surfaces it yields shortest geodesic paths~\cite{kimmel1998triangulated}. In robotics it is used for optimal path planning and navigation in a cost field. For a broader survey of these applications and of the numerical techniques used to solve the equation, we refer the reader to the monograph~\cite{sethian1999book} and to the recent review of Hamilton-Jacobi solvers~\cite{meng2025review}.
 
A large number of mesh-based numerical methods has been developed to find the viscosity solution of~\eqref{eq:eikonal_main}, organized around how each enforces upwinding and causality on a grid. Marching (Dijkstra-like, single-pass) methods order the nodes by increasing solution value and accept them once: the one-pass trajectory algorithm of~\cite{tsitsiklis1995efficient} and the Fast Marching Method (FMM)~\cite{sethian1996marching, sethian1999fast} compute the solution in $O(N \log N)$ for $N$ nodes. PDE-based iterative methods instead sweep the grid in alternating directions with Gauss-Seidel updates until a fixed point is reached: the Fast Sweeping Method (FSM)~\cite{tsai2003sweeping, kao2005sweeping, zhao2005fast} attains $O(N)$ on Cartesian grids using a Godunov numerical Hamiltonian. Higher-order accuracy on these iterative solvers is obtained with WENO reconstructions~\cite{zhang2006high} and with discontinuous-Galerkin sweeping~\cite{li2008dg}. The schemes extend beyond Cartesian grids to general and unstructured meshes, both for marching on triangulated domains~\cite{kimmel1998triangulated, sethian2000unstructured} and for sweeping on triangular meshes~\cite{qian2007triangular, zhang2003weno}, the Ordered Upwind Methods of~\cite{sethian2003ordered} treat anisotropic static Hamilton-Jacobi equations on such meshes. On fully general three-dimensional polyhedral meshes, the finite volume framework of~\cite{hahn2022soner} solves a bidirectional time-relaxed form of the eikonal equation with the Soner boundary condition, while~\cite{hahn2024laplacian} replaces the time relaxation by a Laplacian (vanishing-viscosity) regularization, both recovering the viscosity solution. Semi-Lagrangian schemes discretize the dynamic-programming form by following characteristics back to a previous front and interpolating, with a unified treatment given by~\cite{falcone2014semilagrangian} and the connection to discrete representation formulae and Godunov methods made explicit in~\cite{falcone2002semilagrangian}. Finally, a variational route casts the distance to $\Gamma$, and more generally the solution of the eikonal equation, as the solution of an optimization problem rather than of the equation in differential form. One branch recovers the distance function as the minimizer of an energy, as in the $p$-Poisson formulation that approaches the distance as the order grows~\cite{belyaev2015variational}; in the same spirit, the eikonal equation is solved by elliptic solvers through a penalty and biharmonic regularization~\cite{glowinski2015variational}. Another route casts the problem as a Kantorovich-Rubinstein problem, the maximization of a linear functional over $1$-Lipschitz fields~\cite{belyaev2020admm, ennaji2021augmented}. In high dimensions, where a grid is infeasible, related representation-formula solvers based on the Hopf formula evaluate the solution pointwise as a finite dimensional optimization and thereby sidestep the curse of dimensionality~\cite{darbon2016hopf, chow2017hopf}. What all of these families share, and what matters for the present work, is that each is tied to an explicit spatial discretization of the domain and is solved afresh for every new pair $(\Gamma, F)$.
 
A more recent line of work solves the eikonal equation with scientific machine learning. Most of these methods are meshless: a neural network is trained on the pointwise eikonal residual, or a closely related geometric loss, at sampled collocation points, with spatial derivatives taken by automatic differentiation, so no grid is required~\cite{bin2021pinneik, grubas2023neural, taufik2022upwind, park2024ReSDF}. Meshless viewpoint underlies implicit neural surface reconstruction, where the eikonal equation is solved for the opposite unknown: the velocity is fixed to unit slowness and the source set, the surface, is itself recovered as the zero level set of a learned signed distance function~\cite{park2019deepsdf, atzmon2020sal, gropp2020implicit}, with later work targeting the viscosity solution explicitly~\cite{yang2023steik, krishnan2025viscoreg}. A complementary direction instead brings a mesh back in to attack the curse of dimensionality: the finite difference least-squares method~\cite{EsteveYague2025} enforces a monotone, viscosity-consistent stencil at collocation points sampled on the fly, reaching dimensions where a fixed grid is infeasible. Setting this high-dimensional case aside, the residual-trained meshless solvers share two limitations as forward solvers. First, each model represents the solution of a single instance: a new source set or velocity field requires retraining from scratch, so the amortization that would motivate a learned solver is not realized. Second, the PDE residual does not by itself approximate the viscosity solution: away from the source set many fields satisfy the eikonal equation, and a residual-trained network can converge to a non-physical solution; see the physics-informed neural network (PINN) solvers of the eikonal equation~\cite{grubas2023neural}.

Operator learning is the standard remedy for the first limitation, learning the solution map $\mathcal{S}: (\Gamma, F) \mapsto u$ itself rather than a single solution to use a branch-trunk DeepONet~\cite{lu2021deeponet}, a Fourier Neural Operator (FNO)~\cite{li2021fno}, or their physics-informed, label-free training~\cite{wang2021pideeponet}. Within eikonal and Hamilton-Jacobi problems, the CNN-enhanced DeepONet~\cite{mei2024fully} learns the $\Gamma \mapsto u$ map but is trained supervised on fast-marching targets, inheriting the precomputation cost of the solver it replaces, whereas the physics-informed Fourier neural operator~\cite{song2024seismic} and the Fourier-DeepONet~\cite{zhu2023fourier} train on the residual alone but evaluate it with central differences. The closest antecedent is the Physics-Informed Fourier-DeepONet~\cite{liu2025physics}, which uses the same DeepONet-FNO backbone and residual-only training that we adopt. It encodes a single source through a homogeneous-medium background traveltime $T_0$, which both ties the model to one source per prediction and lets it sidestep the soft boundary penalty that the authors report to be delicate to balance against the residual.

The proposed method, the Lower-envelope Eikonal Neural Operator (LENO), removes both limitations, keeping this operator backbone and physics-informed training, but adds two structural ingredients tailored to the eikonal equation. First, a lower-envelope output transform writes the prediction as the pointwise minimum of per-source candidate fields, each playing the role of the weighted distance from one source. Since the multi-source viscosity solution is itself the pointwise minimum of these distances~\cite{BardiCapuzzoDolcetta1997}, the transform constrains the prediction to the lower-envelope form of the viscosity solution. It also enforces the source condition by construction, so no boundary penalty enters the loss, and reproduces the multi-source structure of the solution within one trained model. Second, the residual in the loss is discretized with the Godunov upwind scheme~\cite{rouy1992viscosity, zhao2005fast, osher1988fronts}, the monotone upwind discretization on which classical eikonal solvers rely to converge to the viscosity solution~\cite{BarlesSouganidis1991}, rather than the central differences used above.

We validate these design choices through experiments on multiple point sources, arbitrary open and closed curve sources, sample-dependent piecewise velocity fields, and realistic OpenFWI CurveVel-A and Style-A media. Additional studies further isolate the effect of the residual discretization, the source encoder, the velocity-augmented trunk input, and the lower-envelope transform, and examine computational efficiency and resolution transfer. These results show that LENO is not merely a residual-trained surrogate for a fixed configuration, but a source- and velocity-conditioned operator whose structural components are essential for heterogeneous multi-source eikonal problems.

\section{Method}\label{sec:method}
In this work, the source set $\Gamma$ is either a finite set of points or an
arbitrary curve (open or closed). It is represented by a finite collection of source points $\{s_k\}_{k=1}^K\subset\Gamma$. Rather than solving a particular instance of~\eqref{eq:eikonal_main}, we learn the solution operator
\begin{equation}
	\label{eq:solution_operator_main}
	\mathcal{S}:(\Gamma,F)\longmapsto u
\end{equation}
across source configurations and velocity fields, so that after training an approximation of the viscosity solution is produced by a single forward pass for any new pair $(\Gamma,F)$.

The design of the operator $\mathcal{S}$ is guided by the variational representation of the viscosity solution $u$ as a weighted distance. Defining the slowness-weighted length of a Lipschitz path $\gamma:[0,1]\to\bar\Omega$ and the induced distance
\begin{align}
    \label{eq:weighted_dist_xy}
    d_F(x,y) = \inf_{\gamma} \left\{L_F(\gamma) \:\middle|\: \gamma(0) = y, \: \gamma(1) = x\right\}, \quad L_F(\gamma) = \int_0^1 \frac{\|\dot\gamma(\tau)\|}{F(\gamma(\tau))}d\tau,
\end{align}
the solution $u$ is the value function of minimization~\cite{BardiCapuzzoDolcetta1997, CrandallLions1983, DeckelnickElliott2004}:
\begin{equation}
	\label{eq:viscosity_representation_main}
	u(x)=\inf_{\gamma} \left\{L_F(\gamma) \:\middle|\: \gamma(0) \in \Gamma, \: \gamma(1) = x\right\} = \min_{y\in\Gamma} d_F(x,y).	
\end{equation}
This representation shows that $u$ is the lower envelope of the source-indexed candidate fields $d_F(\cdot,y)$, $y\in\Gamma$. 

Two structural properties follow directly from this form. First, $u$ is smooth wherever the nearest source is unique and has a kink wherever it is not, so it is only Lipschitz; the resulting nonsmooth set is the cut locus, which reduces to the set of points equidistant (in $d_F$) from two or more nearest sources when $F$ is uniform and $\Gamma$ consists of isolated points. A discontinuous $F$ introduces a physical mechanism for kinks rather than a geometric one: where $F$ jumps from one region to the next, $u$ kinks along that interface. Second, each candidate vanishes at its own source, $d_F(y,y)=0$, so the envelope is zero exactly on $\Gamma$; the boundary condition $u|_\Gamma=0$ is thus intrinsic to the representation rather than a constraint to be imposed. These two properties explain why smooth-prior approaches, such as the standard PINN, are fragile. Such an approach represents $u$ by a single smooth field and imposes $u|_\Gamma=0$ through a soft penalty. A smooth ansatz must then approximate a nonsmooth solution, with the boundary penalty balanced against the interior residual throughout training. One field is asked to do three things at once: form the cone-like singularity at each source, approximate the viscosity solution, and capture the nonsmooth competition among the candidates $d_F(\cdot,y)$, $y\in\Gamma$. Forcing a single smooth network to meet all three demands at once is what makes such an approach challenging to train."

We instead build this lower-envelope structure directly into the operator through three components: i) a DeepONet-FNO backbone (Section~\ref{subsec:deeponet_fno_expanded}) produces, from $\Gamma$ and $F$, one candidate field per source sample $s_k$, playing the role of the weighted distance from that sample; ii) a lower-envelope output transform (Section~\ref{subsec:hardmin_expanded}) takes their pointwise minimum, which enforces $u|_\Gamma=0$ exactly and removes the boundary penalty from the loss; iii) a Godunov upwind residual (Section~\ref{subsec:deeponet_fno_expanded}) provides the sole, label-free training signal. Being a monotone, consistent Godunov discretization of $F\|\nabla u\|$, this residual stays bounded at any finite grid size, even at the nonsmooth cut locus, so the training loss does not blow up. As the grid size tends to zero, the discretization converges to the viscosity solution~\cite{osher1988fronts,rouy1992viscosity}.

\subsection{Representation of source sets and velocity fields}
\label{subsec:representation_expanded}
The source set $\Gamma$ enters the model only through a finite collection of source samples $\{s_k\}_{k=1}^{K}\subset\bar{\Omega}$ obtained by uniform sampling of $\Gamma$. Since the cardinality $K$ may vary across instances of the operator-learning problem, every sample collection is zero-padded to a prescribed maximum cardinality $K_{\text{M}}\ge K$, and a binary mask $m=(m_1,\ldots,m_{K_{\text{M}}})\in\{0,1\}^{K_{\text{M}}}$ is retained alongside it, with $m_k=1$ marking valid samples and $m_k=0$ otherwise.

Rather than feeding the raw coordinates of $\{s_k\}$ to the network, we encode $\Gamma$ as a stack of Euclidean distance fields. For each valid or padded slot $s_k$, define
\begin{equation}
\label{eq:distance_stack_expanded}
    D_k(x)=\|x-s_k\|, \qquad x\in\Omega_h, \qquad k=1,\ldots,K_{\text{M}},
\end{equation}
where $\Omega_h$ is a uniform Cartesian grid of size $N_x\times N_y$, and these fields are stacked into
\begin{equation}
\label{eq:distance_stack_tensor}
    D=(D_1,\ldots,D_{K_{\text{M}}})\in\mathbb{R}^{K_{\text{M}}\times N_x\times N_y}.
\end{equation}
This representation is deliberately redundant: a coordinate pair $s_k$ contains two numbers, whereas $D_k$ describes how that point is seen from every query location of $\Omega_h$. Since the solution itself is a weighted distance from $\Gamma$, supplying per-source distance fields gives the network a geometrically meaningful coordinate system before any training occurs and makes the subsequent output transform natural.

The velocity field $F$ is supplied pointwise to the trunk network. At each query point $x=(x_1,x_2)\in\Omega_h$, the trunk input is
\begin{equation}
\label{eq:trunk_feature_expanded}
    z(x)=\bigl(x_1,x_2,F(x),1/F(x)\bigr)\in\mathbb{R}^4.
\end{equation}
The function values $F$ and $1/F$ are deterministically related, hence informationally redundant, but they are the two forms in which the medium enters the problem: $1/F$ is the local slope and accumulated integrand of the solution~\eqref{eq:weighted_dist_xy}, whereas $F$ is the coefficient of the residual $r=F\|\nabla u\|-1$ that is driven to zero. Since the map $F\mapsto 1/F$ is a reciprocal nonlinearity with derivative $-1/F^{2}$, largest where $F$ is small, a ReLU-type trunk reproduces it only approximately and least accurately there. Presenting both therefore lets the network use whichever form a computation requires, instead of synthesizing a division internally.

\begin{figure}
\centering
\begin{tabular}{c}
\includegraphics[width=\linewidth]{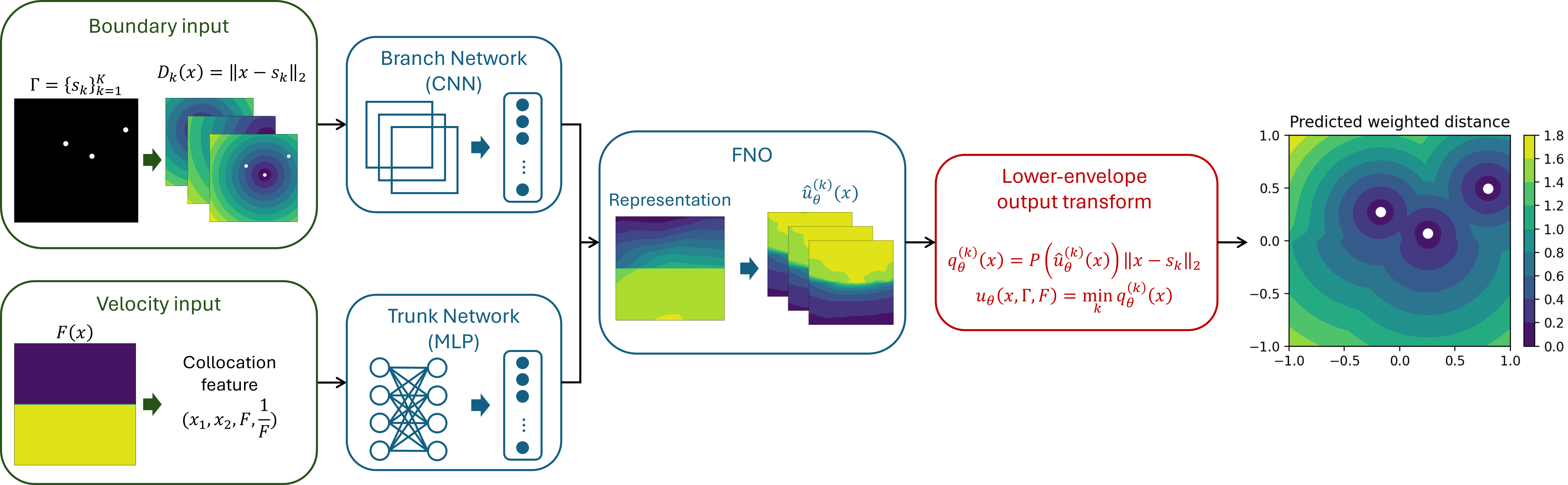}
\end{tabular}
\caption{Overall architecture of LENO. The source set $\Gamma$ is encoded by a CNN-based branch network as a stack of per-source Euclidean distance fields, and the augmented pointwise feature $z(x)=(x_1,x_2,F(x),1/F(x))$ is encoded by an MLP-based trunk network. The branch and trunk representations are coupled by the DeepONet inner product and refined by a Fourier neural operator (FNO) into $K_{\text{M}}$ per-source candidate fields $\{q_{\theta,k}\}$. The prediction $u_\theta$ is obtained by the
lower-envelope output transform, which enforces $u_\theta|_{\Gamma}=0$ by construction.} 
\label{fig:leno}
\end{figure}

\subsection{DeepONet-FNO backbone}
\label{subsec:deeponet_fno_expanded}

LENO comprises three modules, summarized in Figure~\ref{fig:leno}: a convolutional branch encoder that summarizes
the geometry of $\Gamma$, a pointwise trunk network that evaluates
query-dependent features, and an FNO refinement block that converts their
coupling into $K_{\text{M}}$ candidate fields.

\paragraph{Branch encoder.}
The branch network ingests the distance tensor $D$ and encodes the geometry of
$\Gamma$ into a single global latent vector $b_\theta(\Gamma)$. We use stacked
$3\times3$ convolutions with GELU activations and pooling to extract spatial
patterns of the distance stack (the relative arrangement of individual
samples, their separation, the global shape of the sampled curve), followed
by a small fully connected head that compresses these features into a
fixed-length vector. The same vector is reused across all query points of a
given sample and thus acts as a geometry-conditioning code for the
trunk-branch fusion described below.

\paragraph{Trunk network.}
The trunk network is a pointwise multilayer perceptron applied independently
at every grid location $x\in\Omega_h$. It takes the augmented feature $z(x)$
from~\eqref{eq:trunk_feature_expanded} and maps it through ReLU layers to
a latent vector $t_\theta(z(x))$ of the same dimension as $b_\theta(\Gamma)$.
The trunk thus learns a nonlinear feature map on the extended coordinate
space $(x,F(x),1/F(x))$, a representation in which the local velocity and
slowness sit alongside spatial position.

\paragraph{DeepONet coupling.}
The branch and trunk meet through the standard DeepONet inner product,
\begin{equation}
\label{eq:deeponet_inner_expanded}
    h_\theta(x;\Gamma,F)
    = \bigl\langle\, b_\theta(\Gamma),\; t_\theta(z(x))\, \bigr\rangle,
\end{equation}
which couples the global geometric code of $\Gamma$ with the local trunk
feature at $x$ and produces a single scalar field on $\Omega_h$. We read this
field as a $\Gamma$-conditioned latent: it already carries joint information
about the geometry of $\Gamma$ and the velocity, but it collapses the
per-source contributions into a single scalar at each point. The variational
representation in~\eqref{eq:viscosity_representation_main}, by contrast,
prescribes the solution as a minimum over per-source candidates
$d_F(\cdot,y)$, a structure the scalar field $h_\theta$ does not yet expose.
Recovering this per-source separation, and turning the coarse latent into a
stack of nonlinearly refined candidates, is the role of the FNO refinement
introduced next.

\paragraph{FNO refinement and channel lifting.}
The latent $h_\theta$ is passed to a Fourier neural operator module that
first lifts it to multiple channels, applies two spectral convolution blocks
in each spatial direction, and then projects to $K_{\text{M}}$ output channels,
\begin{equation}
\label{eq:fno_candidate_channels_expanded}
    \hat U_\theta(x;\Gamma,F)
    = \bigl(\hat u_{\theta, 1}(x),\ldots,\hat u_{\theta, K_{\text{M}}}(x)\bigr).
\end{equation}
The role of the FNO is twofold. First, its spectral convolutions supply a
global spatial correction to the pointwise DeepONet field, which is necessary
because the candidate that attains the minimum at any location depends on the
global arrangement of all samples and on the velocity field. Second, it lifts
$h_\theta$ from one scalar field into $K_{\text{M}}$ per-source channels. The
$k^{\text{th}}$ channel $\hat u_{\theta,k}$ is not itself the weighted distance
from $s_k$; rather, it is a learned multiplicative correction which, after a softplus and
multiplication by $\|x-s_k\|$ in the lower-envelope output transform of Section~\ref{subsec:hardmin_expanded}, plays the role of $d_F(\cdot,s_k)$.

\subsection{Lower-envelope transform and residual discretization}\label{subsec:hardmin_expanded}
The variational form~\eqref{eq:viscosity_representation_main} writes the
viscosity solution as the lower envelope $u(x)=\min_{y\in\Gamma}d_F(x,y)$,
with $d_F$ from~\eqref{eq:weighted_dist_xy}. We realize this representation in
the architecture through the correspondences
\begin{equation}
\label{eq:correspondence}
\begin{aligned}
    y\in\Gamma \;&\longmapsto\; s_k,\ \ k=1,\ldots,K_{\text{M}},\\
    d_F(\cdot,y) \;&\longmapsto\; q_{\theta,k},\\
    \min_{y\in\Gamma}(\,\cdot\,) \;&\longmapsto\;
        \min_{1\le k\le K_{\text{M}},\,m_k=1}(\,\cdot\,).
\end{aligned}
\end{equation}
The first line is the discretization of $\Gamma$ introduced in
Section~\ref{subsec:representation_expanded}, and the third restricts the
minimum to valid (non-padded) samples through the mask $m$; both are
mechanical. The second line, the design of the per-source candidate
$q_{\theta,k}$, is the centerpiece of the method, and we develop it first.

\paragraph{The per-source candidate.}
Each true candidate $d_F(\cdot,y)$ in~\eqref{eq:viscosity_representation_main}
is a weighted distance from its anchor, and the architecture must be able to
represent it. The admissible form is tightly constrained. Under
$0<F_{\min}\le F\le F_{\max}$ on a convex $\bar{\Omega}$, a straight-segment upper
bound and a length lower bound give, for every $y\in\Gamma$ and $x\neq y$,
\begin{equation}
\label{eq:dist_comparability}
    \frac{1}{F_{\max}}\,\|x-y\|
    \;\le\; d_F(x,y) \;\le\;
    \frac{1}{F_{\min}}\,\|x-y\|,
\end{equation}
so the ratio $d_F(x,y)/\|x-y\|$ lies in $[1/F_{\max},\,1/F_{\min}]$. Every
weighted distance therefore factors as a bounded positive correction times the
Euclidean distance to its anchor, and conversely every member of this factored
class with a positive correction shares the structural features of a weighted
distance: it vanishes at the anchor, grows conically there, and is nonnegative
everywhere. We adopt exactly this class, anchoring at $s_k$ and parameterizing
the correction as the softplus of the FNO output channel,
\begin{equation}
\label{eq:candidate_field_expanded}
    q_{\theta,k}(x)
    = P \left(\hat u_{\theta,k}(x)\right)\,\|x-s_k\|,
    \qquad P(z)=\log(1+e^z)>0.
\end{equation}
The Euclidean factor $\|x-s_k\|$ supplies the conical skeleton that vanishes at
$s_k$, while the learned correction $P(\hat u_{\theta,k})$ sets its scale,
tending to $1/F(s_k)$ as $x\to s_k$ so as to match $\|\nabla u\|=1/F$ at the
anchor and absorbing the spatial variation of $F$ elsewhere. This factorized form is compatible with the structure of a weighted distance: it is nonnegative, vanishes at its anchor, and can represent spatial variations of the medium through the learned positive correction. Hence $q_{\theta,k}$ plays the role of $d_F(\cdot,s_k)$, and their lower
envelope,
\begin{equation}
\label{eq:hardmin_expanded}
    u_\theta(x;\Gamma,F)
    = \min_{1\le k\le K_{\text{M}},\,m_k=1}\, q_{\theta,k}(x),
\end{equation}
realizes $u$ in~\eqref{eq:viscosity_representation_main}.

\paragraph{Properties inherited from the architecture.}
Three structural properties of~\eqref{eq:viscosity_representation_main}
transfer to~\eqref{eq:hardmin_expanded} by construction, prior to any training.
\emph{(i) Exact Dirichlet condition.} For any valid sample $s_j$, the candidate anchored there vanishes,
$q_{\theta,j}(s_j)=0$, while every candidate is nonnegative,
$q_{\theta,k}(s_j)\ge0$ for all $k$; hence the minimum is attained at zero,
\begin{equation}
\label{eq:exact_dirichlet}
    u_\theta(s_j;\Gamma,F)=0,
\end{equation}
for every value of $\theta$, mirroring $d_F(y,y)=0$
in~\eqref{eq:viscosity_representation_main}; the homogeneous Dirichlet boundary condition is encoded in the architecture and need not appear in the loss. Since the boundary condition is satisfied by~\eqref{eq:exact_dirichlet}, the optimizer does not require a separate source penalty, and the classical PINN trade-off between residual and boundary losses~\cite{son2023enhanced,wang2022and} does not arise.
\emph{(ii) Kinks at no extra cost.} The pointwise minimum of $C^1$ candidates
is nondifferentiable exactly where the active index switches with a gradient
jump: the Clarke subdifferential~\cite[Ch.~2]{clarke1990} of $u_\theta$ equals the convex hull of the
gradients of the candidates attaining the minimum, which is nondegenerate iff
those gradients differ. The capacity to form cut-locus kinks is therefore
built into the pointwise minimum and needs no supervision; that the kinks sit
at the true cut locus follows separately, from~(iii) below and from the
uniqueness of the eikonal viscosity solution. \emph{(iii) Local shaping of each
candidate.} Backpropagation through the $\min$ feeds gradients only to the
candidate attaining the minimum, so each $q_{\theta,k}$ is shaped on its own
region of activity (the discrete Voronoi-type cell of $s_k$ in the learned
metric); outside that region it does not attain the minimum, receives no
gradient, and is otherwise unconstrained.

\paragraph{Godunov upwind residual.}
With the architecture supplying~\eqref{eq:exact_dirichlet} for free, the training signal is the eikonal residual evaluated on a uniform grid of spacing
$(\Delta x,\Delta y)$. Let $(i,j)$ index a grid node
$x_{i,j}\in\Omega_h$ with $1\le i\le N_x$ and $1\le j\le N_y$, and write
$u_{i,j}:=u_\theta(x_{i,j};\Gamma,F)$. Denote by $D_x^{\pm}$ and $D_y^{\pm}$
the standard forward ($+$) and backward ($-$) finite difference operators in
$x$- and $y$-directions on $\Omega_h$, with missing differences at boundary nodes replaced
by the nearest valid one. The Godunov approximation to $\|\nabla u\|$ and the
discrete residual at $(i,j)$ are
\begin{align}
    G_{i,j}(u) &= \Bigl[\,
        \max\bigl\{(D_x^{-}u_{i,j})_+,(-D_x^{+}u_{i,j})_+\bigr\}^{2}
        +\max\bigl\{(D_y^{-}u_{i,j})_+,(-D_y^{+}u_{i,j})_+\bigr\}^{2}
        +\varepsilon
    \,\Bigr]^{1/2}, \label{eq:godunov_expanded} \\
    r_{i,j}(u_\theta,F) &= F_{i,j}\,G_{i,j}(u_\theta)-1, \label{eq:godunov_res}
\end{align}
with $(a)_+=\max\{a,0\}$ and $\varepsilon=10^{-12}$. As a monotone, consistent,
one-sided discretization of $F\|\nabla u\|$, this residual is consistent with
the viscosity solution in the limit of the grid size to zero~\cite{osher1988fronts,sethian1999book, rouy1992viscosity}. The Godunov residual is only first-order accurate, but a formally second-order centered difference scheme brings no compensating benefit: lacking monotonicity, it does not converge to the viscosity solution. 

\subsection{Training objective}
\label{subsec:training_expanded}
For a mini-batch $\mathcal{B}$, let $u_\theta^{(n)}=u_\theta(\cdot;\Gamma^{(n)},F^{(n)})$ for $n\in\mathcal{B}$ denote the predictions on $\Omega_h$, and write
\begin{equation}
\label{eq:residual_pointwise}
    r_{i,j}^{(n)}
    = F^{(n)}\bigl((x_{i,j}\bigr)G_{i,j}\bigl(u_\theta^{(n)}\bigr)-1
\end{equation}
for the pointwise eikonal residual at grid node $(i,j) \in \Omega_h$. We train by
minimizing
\begin{equation}
\label{eq:loss_expanded}
    \mathcal{L}(\theta)
    =\frac{1}{|\mathcal{B}||\Omega_h|} \sum_{n\in\mathcal{B}} \sum_{(i,j)\in\Omega_h}
      \left[\left(r_{i,j}^{(n)}\right)^{2}
      +\ell_\delta\!\left(r_{i,j}^{(n)}\right)\right],
\end{equation}
where $\ell_\delta$ is the Huber (smooth-$L^1$) loss. The squared term
penalizes the global residual energy, while the smooth-$L^1$ term provides an
auxiliary residual penalty that is quadratic near zero and grows linearly for
larger residuals. In practice, this mixed objective provides a stable residual
training signal near source samples, cut loci, and velocity interfaces. The loss contains neither a supervised mismatch term nor a boundary penalty: the boundary condition $u|_\Gamma=0$ is satisfied exactly by~\eqref{eq:hardmin_expanded}, and fast-marching
solutions are generated only for post-training evaluation.

\begin{algorithm}
\caption{Physics-informed training of LENO}
\label{alg:self_supervised_training_expanded}
\begin{algorithmic}[1]
\REQUIRE Source/velocity pairs $\{(\Gamma^{(n)},F^{(n)})\}_{n=1}^N$, grid $\Omega_h$,
maximum cardinality $K_{\text{M}}$, batch size $B$, learning rate $\eta$.
\STATE Initialize network parameters $\theta$.
\REPEAT
    \STATE Sample a mini-batch $\mathcal{B}$ of size $B$ uniformly from $\{1,\ldots,N\}$.
    \FOR{$n\in\mathcal{B}$}
        \STATE Load source samples $\Gamma_{n,h}=\{s_1,\ldots,s_{K_n}\}$ and zero-pad to $K_{\text{M}}$ with mask $m^{(n)}$.
        \STATE Construct distance fields $D_k(x)=\|x-s_k\|$ on $\Omega_h$ for $k=1,\ldots,K_{\text{M}}$.
        \STATE Encode $\{D_k\}_{k=1}^{K_{\text{M}}}$ by the branch network to obtain $b_\theta(\Gamma^{(n)})$.
        \STATE Evaluate the trunk features $t_\theta(z(x))$ for all $x\in\Omega_h$, where $z(x)=(x_1,x_2,F^{(n)}(x),1/F^{(n)}(x))$.
        \STATE Form $h_\theta$ via the DeepONet inner product~\eqref{eq:deeponet_inner_expanded} and lift it to correction channels $\{\hat u_{\theta,k}\}_{k=1}^{K_{\text{M}}}$ by the FNO.
        \STATE Form the candidates $q_{\theta,k}$ via~\eqref{eq:candidate_field_expanded} and apply the lower-envelope transform~\eqref{eq:hardmin_expanded} to obtain $u_\theta^{(n)}$.
        \STATE Compute the Godunov residual $r_{i,j}^{(n)}$ on $\Omega_h$ via~\eqref{eq:residual_pointwise}.
    \ENDFOR
    \STATE Assemble $\mathcal{L}(\theta)$ from $\{r_{i,j}^{(n)}\}_{n\in\mathcal{B}}$ via~\eqref{eq:loss_expanded}.
    \STATE Update $\theta$ by one Adam step with learning rate $\eta$.
\UNTIL{convergence}
\end{algorithmic}
\end{algorithm}

All models are trained with Adam at learning rate $10^{-5}$ and batch size
$200$. The same architectural components are used across experiments; the
only quantity that varies with the source class is the maximum cardinality
$K_{\text{M}}$. The full procedure is summarized in
Algorithm~\ref{alg:self_supervised_training_expanded}.

\section{Numerical experiments}
\label{sec:experiments}
We validate the proposed framework on three groups of eikonal equations. The first two groups are controlled synthetic studies that vary the source geometry, isolated point sources and arbitrary curves (open or closed), under two velocity regimes: a fixed smooth velocity field (E1) and a sample-dependent piecewise constant velocity field (E2). Together, E1 and E2 probe the two axes of the
solution operator $(\Gamma,F)\mapsto u$. E1 isolates the effect of source geometry under a fixed medium, whereas E2 additionally requires conditioning on a function-valued, non-smooth velocity input. The third group (E3) moves beyond synthetic data to realistic velocity models  drawn from the CurveVel-A and Style-A families of the OpenFWI benchmark~\cite{deng2022openfwi}, under a surface-source acquisition geometry of the kind encountered in seismic travel-time computation.

All experiments are conducted on the square domain $\Omega=[-1,1]^2$. The synthetic groups E1 and E2 use a uniform $100\times100$ grid, while the OpenFWI experiments inherit the native $70\times70$ resolution of
the corresponding OpenFWI velocity fields. For each setting, we generate $5{,}000$ configurations, of which $4{,}000$ are used for training and $1{,}000$ are reserved for testing. Reference solutions are computed on the same grid as the corresponding experiment using the \texttt{scikit-fmm} implementation of the fast marching method~\cite{sethian1999fast}. These reference solutions are used only for post-training evaluation and visualization, and never used in the training objective, which is based solely on the Godunov residual.
\subsection{Experimental settings}
\label{subsec:experimental_settings}
\paragraph{E1: Fixed smooth velocity field.}
In the first group, the velocity field is fixed across all samples and
is given by
\begin{equation}
\label{eq:smooth_velocity}
    F(x_1,x_2)=2x_1^2+1,
\end{equation}
which is smooth, strictly positive, and spatially heterogeneous. Since the velocity depends only on the horizontal coordinate, the corresponding viscosity solution contours are distorted nonuniformly across the domain. Within this regime we consider two source geometries. In the point-source case, $\Gamma$ consists of three points sampled in $[-1,1]^2$, so that the maximum cardinality is $K_{\text{M}}=3$. In the curve-source setting, $\Gamma$ is an arbitrary curve, which may be open or closed. Since each curve is represented by a discrete set of sampled source points, the number of source samples may vary across configurations. We therefore pad the sampled source points to a fixed maximum cardinality, set to $K_{\text{M}}=100$ in our experiments. Both cases use the same distance-stack representation for $\Gamma$ and the same neural operator architecture, the only change being the padding size $K_{\text{M}}$. This group serves as a baseline test for two essential properties of the proposed method: its ability to represent the competition among multiple expanding fronts, and its geometric extensibility from a small finite source set to a dense discrete approximation of a continuous source.

\paragraph{E2: Sample-dependent piecewise velocity field.}
The second group retains the two source geometries of E1, namely three point sources ($K_{\text{M}}=3$) and arbitrary curves ($K_{\text{M}}=100$), but replaces the fixed smooth velocity by a sample-dependent piecewise constant velocity field. In contrast to E1, the model encodes the heterogeneous velocity field \(F\) and the source geometry. The corresponding viscosity solution may exhibit abrupt changes in its gradient across velocity interfaces, and the travel-time contours must bend according to the local propagation speed. This group therefore evaluates whether the trunk features $(x_1, x_2, F(x), 1/F(x))$, together with the FNO refinement module, can learn a genuine operator dependence on the heterogeneous velocity field, and whether this dependence remains effective when the source geometry
varies simultaneously with the medium.

\paragraph{E3: OpenFWI benchmark with surface sources.}
The third group evaluates the framework on velocity models drawn from 
the CurveVel-A and Style-A families of the OpenFWI benchmark~\cite{deng2022openfwi}, 
mixed in equal proportion. Each velocity model is used at its native resolution, 
mapped onto the computational domain $\Omega=[-1,1]^2$. For each of the $5{,}000$ 
configurations, a single source point is placed at a random horizontal position 
along the top boundary $\{x_2=1\}$, mimicking the surface acquisition configuration 
common in seismic travel-time applications.

This group is designed as a direct counterpart to~\cite{liu2025physics}, which 
adopts the same surface-source setting but relies on an imprecise assumption $\tau_{\text{surface}}=0$. This  assumption requires the medium to be approximately homogeneous near the surface, 
which may not hold in realistic models. The lower-envelope transform of LENO  enforces $u|_\Gamma=0$ exactly by construction, without any homogeneous-medium  assumption, and is therefore applicable to arbitrarily heterogeneous velocity  fields.

\paragraph{Evaluation protocol.}
We evaluate the trained models using two complementary criteria. The
first is the test residual, computed with the same Godunov
discretization used during training:
\begin{equation}
\label{eq:test_residual}
    \mathcal{R}_{\mathrm{test}}
    =
    \mathbb{E}_{(\Gamma,F)}
    \left[
    \frac{1}{|\Omega_h|}
    \sum_{(i,j)\in\Omega_h}
    \left\{
        \bigl(r_{i,j}\bigr)^2
        +
        \ell_\delta\!\left(r_{i,j}\right)
    \right\}
    \right],
\end{equation}
where $r_{i,j}=F(x_{i,j})G_{i,j}(u_\theta)-1$ denotes the discrete eikonal residual and $\ell_\delta$ is the smooth-$L^1$ loss. The residual quantifies satisfaction of the discrete eikonal equation.
The second criterion is the relative $L^2$ error against the reference solutions of FMM:
\begin{equation}
\label{eq:fmm_rel_error}
    \mathcal{E}_{\mathrm{rel}}
    =
    \mathbb{E}_{(\Gamma,F)}
    \left[
    \frac{
        \left(\sum_{(i,j)\in\Omega_h}
            \left|u_\theta(x_{i,j})-u_{\mathrm{FMM}}(x_{i,j})\right|^2
        \right)^{1/2}
    }{
        \left(\sum_{(i,j)\in\Omega_h}
            \left|u_{\mathrm{FMM}}(x_{i,j})\right|^2
        \right)^{1/2}
    }
    \right].
\end{equation}
The relative error normalizes by the magnitude of the reference solution, so that errors are comparable across configurations with different velocity scales and source geometries. Since the solutions of FMM are not used during training, this metric serves solely as an external validation criterion.

\begin{table}[t]
\centering
\caption{Evaluation of LENO across the three experiment groups. The model is trained only
with the Godunov residual, and \texttt{scikit-fmm} solutions are used
exclusively for test-set evaluation. Groups E1 and E2 each span two
source geometries (point sources and arbitrary curves) under a fixed
smooth and a sample-dependent piecewise velocity field, respectively;
group E3 uses realistic OpenFWI CurveVel-A and Style-A velocity models with
surface sources.}
\label{tab:expanded_main_results}
\begin{tabular}{llllcc}
\toprule
& Setting & Source & Velocity & $\mathcal{R}_{\mathrm{test}}$
& $\mathcal{E}_{\mathrm{rel}}$ \\
\midrule
\multirow{2}{*}{E1}
   & Point sources & 3 pts & fixed $F=2x_1^2+1$
   & $2.12 \cdot 10^{-3}$ & $2.22 \cdot 10^{-2}$ \\
   & Curve source  & curve ($K_{\text{M}}{=}100$) & fixed $F=2x_1^2+1$
   & $1.07 \cdot 10^{-2}$ & $2.10 \cdot 10^{-2}$ \\
\midrule
\multirow{2}{*}{E2}
   & Point sources & 3 pts & piecewise, per-sample
   & $5.26 \cdot 10^{-3}$ & $3.64 \cdot 10^{-2}$ \\
   & Curve source  & curve ($K_{\text{M}}{=}100$) & piecewise, per-sample
   & $1.17 \cdot 10^{-2}$ & $2.28 \cdot 10^{-2}$ \\
\midrule
E3 & Surface sources & 1 pt on $\{x_2{=}1\}$ & CurveVel-A and Style-A
   & $5.24 \cdot 10^{-3}$ & $4.33 \cdot 10^{-2}$ \\
\bottomrule
\end{tabular}
\end{table}

\subsection{E1: Fixed smooth velocity field}
\label{subsec:e1_results}
\begin{figure}
\centering
\begin{tabular}{c}
\includegraphics[width=\linewidth]{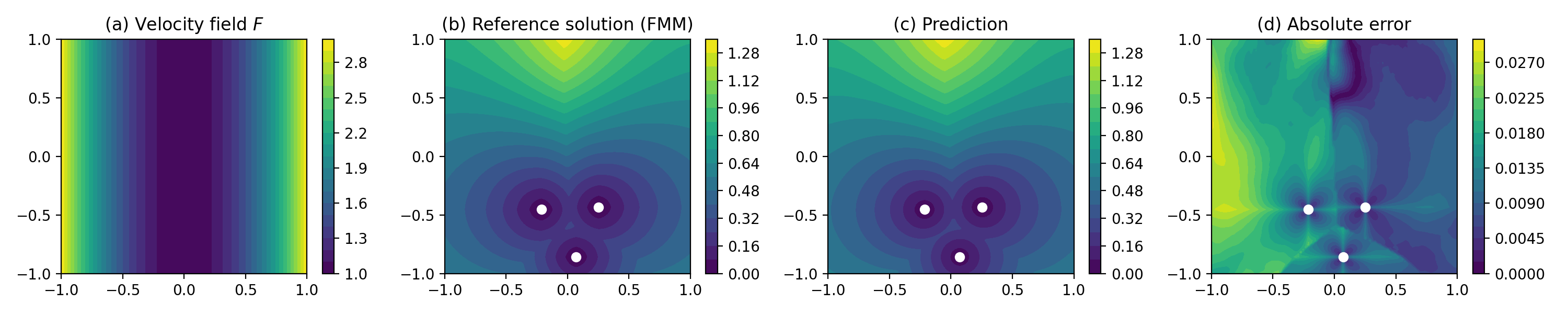} \\
\includegraphics[width=\linewidth]{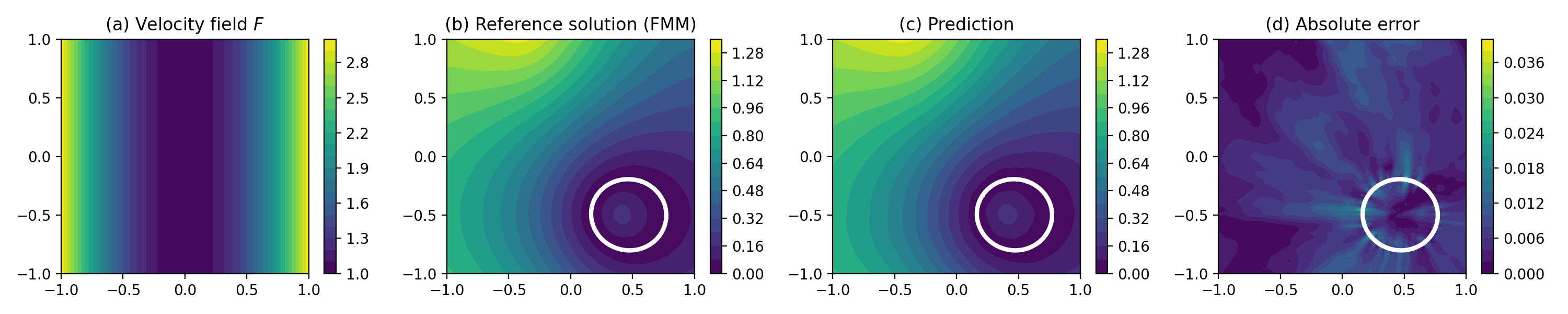}
\end{tabular}
\caption{Qualitative results for the fixed smooth velocity setting E1, where $F(x_1,x_2)=2x_1^2+1$. Top row: three point sources. 
Bottom row: arbitrary curve source. In each row, panels (a)-(d) show, respectively, (a) the velocity field $F$, (b) the \texttt{scikit-fmm} reference solution, (c) the predicted weighted distance field $u_\theta$, and (d) the absolute error. The source set is marked in white in panels (b)-(d), either as points or as a curve. The prediction satisfies the boundary condition by construction through the hard lower-envelope output transform.}
\label{fig:e1_results}
\end{figure}

The top row of Figure~\ref{fig:e1_results} presents a representative test
sample from the point-source case of E1. The predicted travel-time field
vanishes at the three source points and increases away from them, as required
by the Dirichlet condition. The level sets display the expected multi-front
structure: near each source, the contours are organized around the
corresponding source, whereas transition regions form where the active
minimizing source changes. This behavior demonstrates the role of the
lower-envelope output transform, which represents the viscosity solution as
the pointwise minimum of source-wise candidate fields.

The bottom row of Figure~\ref{fig:e1_results} shows a representative test
sample for the curve-source case under the same velocity field. Here the
source set is an extended one-dimensional object rather than a finite
collection of isolated points. The curve is represented by sampled source
points, and the same hard lower-envelope transform is applied over the valid
sampled sources. The predicted field preserves the zero-travel-time structure
along the curve and produces level sets that propagate outward from
it. This shows that the architecture is not tied to isolated point sources:
the same distance-stack branch representation and lower-envelope output
transform extend naturally from a small finite source set to a dense discrete
approximation of a continuous source.

\subsection{E2: Sample-dependent piecewise velocity field}
\label{subsec:e2_results}
\begin{figure}
\centering
\begin{tabular}{c}
\includegraphics[width=\linewidth]{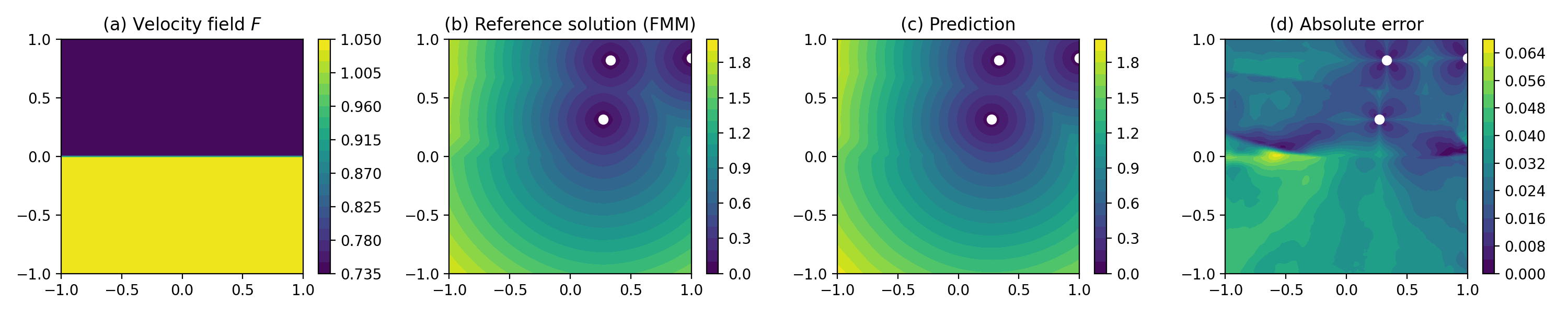} \\
\includegraphics[width=\linewidth]{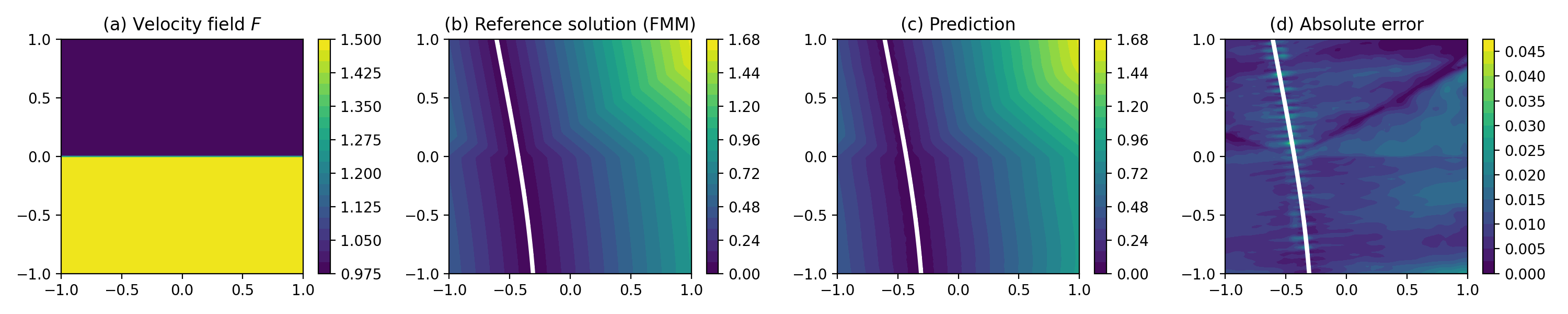} 
\end{tabular}
\caption{Qualitative results for the sample-dependent piecewise constant velocity setting E2. Top row: three point sources. Bottom row: arbitrary curve source. In each row, panels (a)-(d) show, respectively, (a) the velocity field $F$, (b) the \texttt{scikit-fmm} reference solution, (c) the predicted weighted distance field $u_\theta$, and (d) the absolute error. The source set is marked in white in panels (b)-(d), either as points or as a curve. The prediction captures the deformation of travel-time contours induced by the piecewise constant velocity field while satisfying the boundary condition by construction. }
\label{fig:e2_results}
\end{figure}

The top row of Figure~\ref{fig:e2_results} reports a representative test sample
from the point-source case of E2. This setting is more challenging
than E1 because the velocity field varies from sample to sample and is
piecewise constant, so the model must adapt each source configuration to a different propagation medium. In regions of higher velocity, travel
time accumulates more slowly and the level sets are more widely
spaced; in slower regions, the level sets become denser. Across
velocity interfaces, the viscosity solution may exhibit abrupt changes
in its gradient. The predicted field captures the main deformation of
the travel-time contours induced by the heterogeneous velocity, which
indicates that the model is not merely learning a source-to-distance
map for a fixed background medium. Rather, the pointwise velocity
features supplied to the trunk network, together with the global
spatial refinement provided by the FNO, enable the prediction to
depend on the input velocity function. The simultaneous use of $F$ and
$1/F$ is natural in this setting, since the eikonal equation relates
the local gradient magnitude to the reciprocal of the velocity, while
the travel time itself accumulates slowness along paths.

The bottom row of Figure~\ref{fig:e2_results} shows the curve-source case under the
same piecewise velocity regime, in which both the source geometry and
the velocity field vary across samples. The predicted field vanishes
along the sampled curve and propagates outward through the
heterogeneous medium, with level sets that bend across velocity
interfaces in agreement with the fast-marching reference. This is the
most demanding of the synthetic settings, since no single axis of the
operator $(\Gamma,F)\mapsto u$ is held fixed: the model must combine
the curve-source representation validated in E1 with the
velocity-conditioning mechanism validated by the point-source case of
E2. The fact that the same architecture recovers the qualitative structure of the
reference solution, without setting-specific modifications beyond the choice
$K_{\text{M}}=100$, indicates that the geometric and physical components of the
architecture play complementary roles: the distance-stack branch and the
lower-envelope output transform encode the source geometry, while the trunk
velocity features and the FNO refinement account for the heterogeneous medium.

This setting also highlights why an upwind residual is appropriate: the
piecewise constant velocity can produce kinks across interfaces, and the
Godunov discretization is designed to select the one-sided derivatives
consistent with viscosity solutions.

\subsection{E3: OpenFWI benchmark with surface sources}
\label{subsec:e3_results}
\begin{figure}
\centering
 \begin{tabular}{c}
\includegraphics[width=\linewidth]{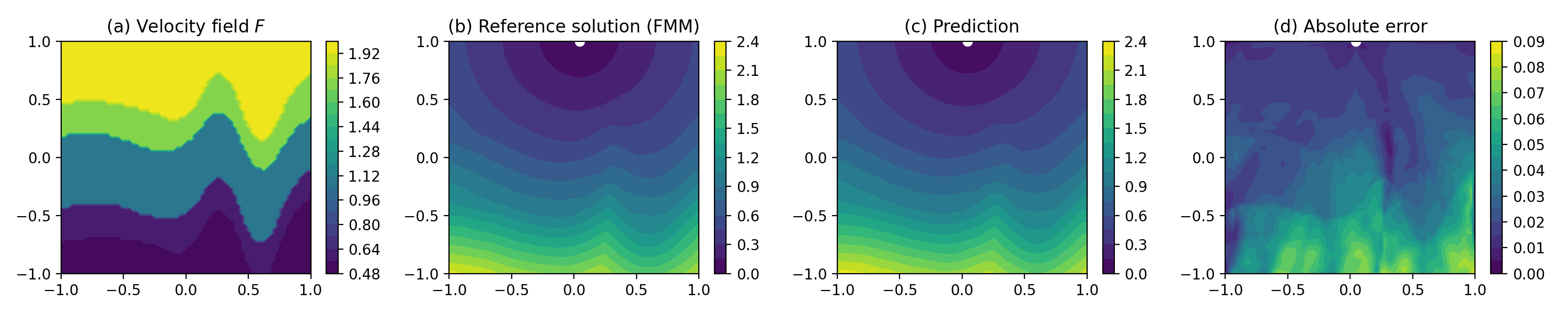} \\
\includegraphics[width=\linewidth]{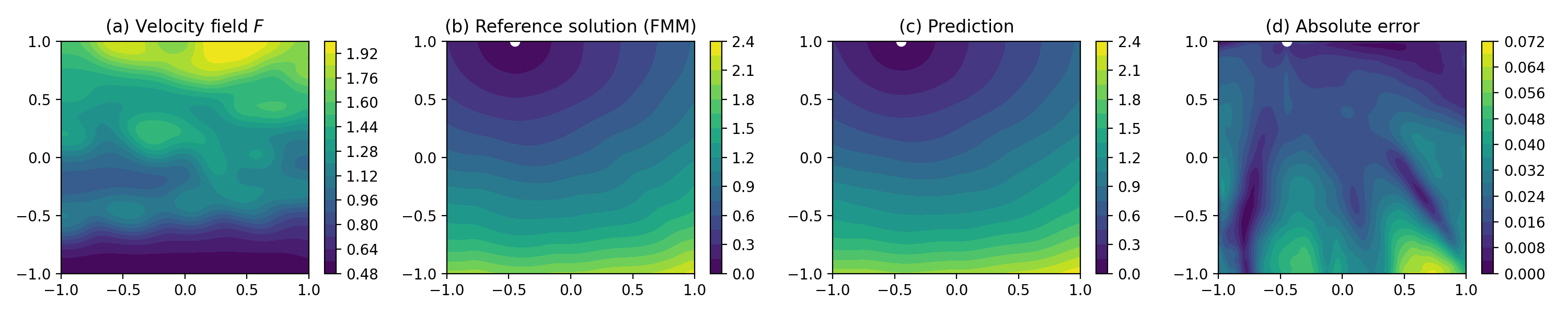}
\end{tabular}
\caption{Qualitative results for the OpenFWI setting E3 with a single random surface source on $\{x_2=1\}$. Top row: a representative sample from the CurveVel-A family. Bottom row: a representative sample from the Style-A family. In each row, panels (a)--(d) show, respectively, (a) the velocity field $F$, (b) the \texttt{scikit-fmm} reference solution, (c) the predicted travel-time field $u_\theta$, and (d) the absolute error. The source point is marked in white in panels (b)-(d). The predicted front propagates downward from the surface source and follows the structures induced by the heterogeneous velocity field, while satisfying the boundary condition exactly by construction, without requiring a homogeneous-medium assumption. As in E1 and E2, the \texttt{scikit-fmm} reference solution is used only for evaluation and visualization, never for training.}
\label{fig:e3_result}
\end{figure}

Figure~\ref{fig:e3_result} presents representative test samples from E3,
one drawn from the CurveVel-A family and one from the Style-A family.
In each case, a single source point is placed on the top boundary
$\{x_2=1\}$ and the predicted travel-time field propagates downward
through the heterogeneous medium. The level sets refract at velocity
interfaces in agreement with the \texttt{scikit-fmm} reference, and
the boundary condition $u=0$ is satisfied exactly at the surface source
by construction.

The key distinction from~\cite{liu2025physics} lies in how the source
condition is enforced. That work decomposes the traveltime as $T = T_0 + \tau$,
where $T_0$ is the background traveltime computed under a homogeneous velocity
model with the velocity value sampled at the source location. The surface
condition is then imposed by setting $\tau_{\text{surface}} = 0$, under the
assumption that the horizontal velocity variation on the surface is negligible
relative to the vertical variation, so that $T_{\text{surface}} \approx
T_{0,\text{surface}}$. The authors themselves acknowledge that this is an
imprecise assumption, which may yield unstable results under complex surface
velocity conditions. LENO requires no such decomposition or assumption: the
lower-envelope transform enforces $u_\theta|_\Gamma = 0$ exactly by
construction, for any velocity field.

Taken together with E1 and E2, this result shows that the proposed operator
handles the same surface-source, realistic-velocity setting as~\cite{liu2025physics}
without the homogeneous-medium prerequisite, while the lower-envelope
architecture provides exact boundary enforcement as a structural property
rather than as a modeling assumption.

\section{Additional validation}
\label{sec:additional_validation}
Beyond the qualitative experiments in Section~\ref{sec:experiments}, we perform four additional studies to further validate the proposed framework.  First, we compare residual discretizations on a high-resolution E1 three-point-source problem to justify the use of the Godunov upwind residual. Second, we conduct ablation experiments on the OpenFWI setting to isolate the roles of the source encoder, velocity-augmented trunk input, and lower-envelope output transform. Third, we compare the inference time of LENO with the fast marching method and the fast sweeping method. Fourth, we examine resolution transfer by training the model on a downsampled $50\times50$ grid and evaluating it on the original $100\times100$ grid. All experiments were conducted on a server equipped with two Intel Xeon Gold 6530 CPUs, 503 GiB of system memory, and NVIDIA RTX 6000 Ada Generation GPUs with 48 GB of memory. Unless otherwise stated, GPU-based experiments used a single NVIDIA RTX 6000 Ada Generation GPU. 

\subsection{High-resolution residual-discretization study}
\label{subsec:stencil_comparison}

We first examine how the choice of finite difference stencil in the physics-informed residual affects the learned eikonal solution. This study is conducted on the E1 three-point-source setting with the fixed smooth velocity
$F(x_1,x_2)=2x_1^2+1$, and the trained models are evaluated on a
high-resolution $1000\times1000$ grid.

We compare three discretizations of the eikonal residual: a forward difference residual, a central-difference residual, and the proposed Godunov upwind residual. Although all three schemes approximate the same formal quantity $F\|\nabla u\|-1$ in smooth regions, the multi-source eikonal solution is only Lipschitz and contains nonsmooth transition sets where the active source changes. Central differences average across such transition sets, while a purely forward stencil is not generally aligned with the causal direction of front propagation. The Godunov residual instead selects the one-sided derivatives compatible with the viscosity solution. Table~\ref{tab:stencil_comparison} reports the test residual and the relative $L^2$ error against the fast-marching reference. The result supports the use of the Godunov upwind residual as the training signal in LENO.

\begin{table}[t]
\centering
\caption{
High-resolution residual-discretization study on the E1 three-point-source setting. All results are evaluated on a $1000\times1000$ grid with the fixed
smooth velocity $F(x_1,x_2)=2x_1^2+1$. The comparison isolates the effect of the finite difference stencil used in the physics-informed eikonal residual, rather than the effect of architectural components.}
\label{tab:stencil_comparison}
\begin{tabular}{lcc}
\toprule
Residual discretization
& $\mathcal{R}_{\mathrm{test}}$
& $\mathcal{E}_{\mathrm{rel}}$ \\
\midrule
Forward difference & $8.83 \cdot 10^{-3}$ & $3.83 \cdot 10^{-2}$ \\
Central difference & $1.57 \cdot 10^{-2}$ & $4.44\cdot 10^{-2}$ \\
Godunov upwind & $5.46 \cdot 10^{-3}$ & \bm{$2.71 \cdot 10^{-2}$} \\
\bottomrule
\end{tabular}
\end{table}

\subsection{Ablation study}
\label{subsec:ablation}
\begin{figure}
\centering
\begin{tabular}{c}
\includegraphics[width=\linewidth]{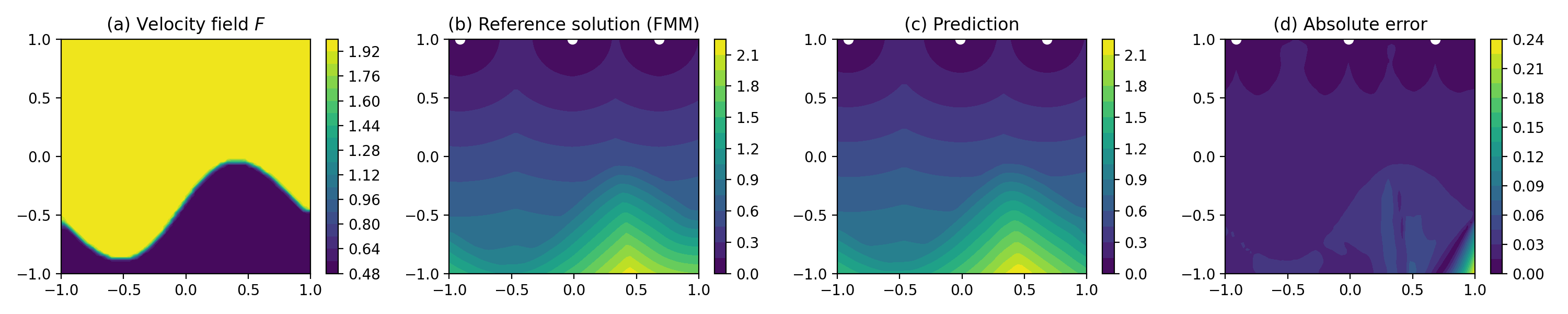}
\end{tabular}
\caption{Representative sample from the FlatFault-A ablation setting with three surface point sources on $\{x_2=1\}$. Panels (a)-(d) show, respectively, (a) the velocity field $F$, (b) the \texttt{scikit-fmm} reference solution, (c) the predicted travel-time field $u_\theta$ of the full LENO model, and
(d) the absolute error. The source points are marked in white in panels (b)-(d). }
\label{fig:ablation_sample}
\end{figure}

Having justified the Godunov upwind residual as the training signal, we next conduct ablation experiments on the OpenFWI FlatFault-A velocity with three surface point sources. This setting provides a stringent test bed because the model must account for both the source geometry and the complex layered, faulted velocity structure. The goal is to isolate the contribution of three key modeling choices: the CNN-based distance-stack branch encoder, the velocity-augmented trunk input, and the lower-envelope output transform. Figure~\ref{fig:ablation_sample} shows a representative prediction of the 
full LENO model on the FlatFault-A setting.

We compare the full model with the following three variants:
\begin{itemize}
    \item \emph{DeepSets branch encoder}: the CNN-based branch encoder that ingests the distance-stack tensor $D=(D_1,\ldots,D_{K_{\text{M}}})$ is replaced by a permutation-invariant DeepSets encoder applied directly to the source point cloud $\{s_k\}_{k=1}^{K}$~\cite{zaheer2017deep}. This variant tests whether the spatial distance-stack representation and convolutional aggregation provide an advantage over a set-based encoder for representing the source geometry. The lower-envelope output transform is kept unchanged.

    \item \emph{w/o trunk feature}: the trunk input is reduced to the spatial
    coordinate $(x,y)$ only, removing the velocity-dependent features $F(x)$
    and $1/F(x)$. This tests whether the model genuinely uses the heterogeneous
    velocity field in the OpenFWI setting, rather than relying only on source
    geometry and spatial coordinates.

    \item \emph{w/o lower-envelope output transform}:  the lower-envelope
    construction is removed, and the model predicts a single travel-time field.
    The boundary condition $u=0$ on $\Gamma$ is then imposed through a soft penalty.
\end{itemize}

\begin{table}[h]
\centering
\caption{
Ablation study on the OpenFWI FlatFault-A setting with three surface point sources. The full model uses a CNN-based branch encoder applied to the
distance-stack source representation, the velocity-augmented trunk input
$(x_1,x_2,F(x),1/F(x))$, and the lower-envelope output transform. The
DeepSets variant replaces only the CNN-based distance-stack branch encoder
with a permutation-invariant point-cloud encoder. In the variant without the
lower-envelope transform, the boundary condition is imposed through a soft
penalty.
}
\label{tab:ablation}
\begin{tabular}{lcc}
\toprule
Model variant & $\mathcal{R}_{\mathrm{test}}$ & $\mathcal{E}_{\mathrm{rel}}$ \\
\midrule
Full LENO & $4.06 \cdot 10^{-3}$ & \bm{$5.80 \cdot 10^{-2}$} \\
DeepSets branch encoder & $2.62 \cdot 10^{-3}$ & $5.91 \cdot 10^{-2}$ \\
w/o trunk feature ($(x,y)$ only) & $2.45 \cdot 10^{-1}$ & $3.38 \cdot 10^{-1}$ \\
w/o lower-envelope output transform & $9.27 \cdot 10^{-3}$ & $3.24 \cdot 10^{-1}$ \\
\bottomrule
\end{tabular}
\end{table}

The ablation results in Table~\ref{tab:ablation} show that the choice of
source encoder has a moderate effect, whereas the velocity-dependent trunk
features and the lower-envelope transform are essential. Replacing the
CNN-based distance-stack branch encoder with a DeepSets point-cloud encoder
yields a comparable test residual, $\mathcal{R}_{\mathrm{test}}=2.62 \cdot 10^{-3}$, but
a slightly larger relative error, $\mathcal{E}_{\mathrm{rel}}=5.91 \cdot 10^{-2}$, compared
with $5.80 \cdot 10^{-2}$ for the full LENO. This suggests that a permutation-invariant
set encoder can represent the three surface sources reasonably well in the
OpenFWI setting, but the distance-stack CNN branch provides a slightly more
accurate geometry-conditioned representation for the final travel-time field.
In contrast, removing the velocity-dependent trunk features causes the
relative error to increase to $3.38 \cdot 10^{-1}$, showing that the model cannot handle the layered and faulted velocity fields using spatial coordinates alone. Removing the lower-envelope output transform similarly degrades the relative error to $3.24 \cdot 10^{-1}$. Since this variant enforces $u=0$ on $\Gamma$ only through a soft boundary loss, the result confirms that the lower-envelope transform is crucial not only for exact source enforcement but also for representing the multi-source viscosity structure of the eikonal solution.

\subsection{Computational efficiency}
\label{subsec:computational_efficiency}

We next compare the inference time of LENO with the fast marching method
(FMM) and the fast sweeping method (FSM). The purpose of this comparison is not to claim that one method is universally superior in all regimes, but to clarify the different computational regimes of marching, sweeping, and learned-operator evaluation. FMM provides a reliable deterministic reference solution for each given source--velocity configuration. However, because it solves each new instance independently, its total runtime grows almost linearly with the number of test instances, as shown in the first row of Table~\ref{tab:runtime}. FSM with CPU reduces this cost by using a simpler sweeping procedure that iteratively updates the solution from neighboring grid values. As shown in Table~\ref{tab:runtime}, FSM with CPU is faster than FMM for both one-sample and 1,000-sample evaluations. Nevertheless, it still evaluates the instances sequentially on the CPU, so its total runtime continues to accumulate as the number of samples increases. FSM with GPU addresses this limitation by applying the sweeping step to many independent source--velocity instances simultaneously. Following the parallel FSM strategy of~\cite{detrixhe2013paralleFSM}, for each sweep direction, grid points on the same diagonal level can be updated in parallel. Our batched implementation further applies these updates across all samples simultaneously. This can lead to a substantial reduction in total runtime compared to both FMM and FSM with CPU, as shown in the third row of Table~\ref{tab:runtime}. 

Finally, LENO does not solve the discrete eikonal problem iteratively at inference time like FMM and FSM. Instead, it evaluates the learned solution operator through a neural-network forward pass, which consists of a fixed sequence of linear and nonlinear mappings and is highly efficient for batched GPU computation. Consequently, the 1,000-sample runtime remains nearly the same as the one-sample runtime in Table~\ref{tab:runtime}. As shown in the last row of Table~\ref{tab:runtime}, for a single sample, LENO is slightly faster than FMM and comparable to FSM with CPU in both the three-point-source (E1) and curve-source cases (E1). However, for 1,000 samples, LENO processes the test set in a batched forward pass, so the total evaluation time remains nearly unchanged at the scale of 1,000 instances. Compared with FMM and the CPU implementation of FSM, LENO reduces the total evaluation time for 1,000 test instances by several orders of magnitude. It also remains substantially faster than the GPU implementation of FSM. This indicates that the main computational advantage of LENO does not come from a large per-instance speed-up, but from batched evaluation over many test instances, and that the trained operator is most advantageous in repeated-query regimes, where solutions must be evaluated for many source configurations or velocity fields.

\begin{table}[h]
\centering
\caption{Computation time comparison among FMM, FSM with CPU, FSM with GPU, and LENO.
All times are reported in seconds (s). The ``1 sample'' columns report the mean $\pm$
standard deviation over 1,000 repeated single-instance evaluations. The ``1,000 samples''
columns report the mean $\pm$ standard deviation over five repeated evaluations of all
1,000 test instances. FMM and CPU-FSM evaluate the test instances sequentially on the CPU,
while GPU-FSM and LENO evaluate them on a single GPU using batched computation.}
\label{tab:runtime}
\begin{tabular}{l *{4}{R@{}L@{}P@{}R@{}L}}
\toprule
\multirow{2}{*}{Method}
  & \multicolumn{10}{c}{Three-point sources}
  & \multicolumn{10}{c}{Curve source} \\
\cmidrule(lr){2-11}\cmidrule(lr){12-21}
  & \multicolumn{5}{c}{1 sample} & \multicolumn{5}{c}{1,000 samples}
  & \multicolumn{5}{c}{1 sample} & \multicolumn{5}{c}{1,000 samples} \\
\cmidrule(lr){1-1}\cmidrule(lr){2-6}\cmidrule(lr){7-11}\cmidrule(lr){12-16}\cmidrule(lr){17-21}
FMM
  & 5.02 & {}\cdot 10^{-3} & {}\pm{} & 3.22 & {}\cdot 10^{-4}
  & 4.99 & {}              & {}\pm{} & 3.25 & {}\cdot 10^{-4}
  & 4.97 & {}\cdot 10^{-3} & {}\pm{} & 3.30 & {}\cdot 10^{-5}
  & 5.00 & {}              & {}\pm{} & 1.93 & {}\cdot 10^{-2} \\
FSM with CPU
  & 1.01 & {}\cdot 10^{-3} & {}\pm{} & 6.80 & {}\cdot 10^{-5}
  & 9.28 & {}\cdot 10^{-1} & {}\pm{} & 8.94 & {}\cdot 10^{-4}
  & \bm{6.82} & \bm{{}\cdot 10^{-4}} & \bm{{}\pm{}} & \bm{1.10} & \bm{{}\cdot 10^{-5}}
  & 8.84 & {}\cdot 10^{-1} & {}\pm{} & 3.73 & {}\cdot 10^{-3} \\
FSM with GPU
  & 2.26 & {}\cdot 10^{-2} & {}\pm{} & 7.52 & {}\cdot 10^{-4}
  & 9.34 & {}\cdot 10^{-2} & {}\pm{} & 2.78 & {}\cdot 10^{-2}
  & 1.72 & {}\cdot 10^{-2} & {}\pm{} & 6.82 & {}\cdot 10^{-4}
  & 7.59 & {}\cdot 10^{-2} & {}\pm{} & 4.29 & {}\cdot 10^{-3} \\
LENO
  & \bm{9.26} & \bm{{}\cdot 10^{-4}} & \bm{{}\pm{}} & \bm{8.81} & \bm{{}\cdot 10^{-5}}
  & \bm{1.03} & \bm{{}\cdot 10^{-3}} & \bm{{}\pm{}} & \bm{9.97} & \bm{{}\cdot 10^{-6}}
  & 9.17 & {}\cdot 10^{-4} & {}\pm{} & 3.25 & {}\cdot 10^{-5} 
  & \bm{1.06} & \bm{{}\cdot 10^{-3}} & \bm{{}\pm{}} & \bm{2.05} & \bm{{}\cdot 10^{-5}} \\
\bottomrule
\end{tabular}
\end{table}

We finally examine whether the learned operator can be evaluated on a finer
output grid than the one used during training. In this experiment, the model is trained on a $50\times50$ grid and then evaluated on the original $100\times100$ grid. The branch CNN input is kept at $50\times50$ resolution, because the branch encoder takes the source-distance stack as a fixed-resolution image. In contrast, the trunk network, the FNO refinement block, and the lower-envelope output transform are applied on the $100\times100$ evaluation grid without architectural modification.

This resolution transfer is possible because the trunk network is evaluated
pointwise at query features $z(x)=(x_1,x_2,F(x),1/F(x))$ and is therefore not
tied to the training grid. The lower-envelope output transform is also pointwise once the source samples and candidate fields are given. The FNO
block is applied to the grid-level latent representation at the evaluation
resolution and produces candidate fields on the corresponding grid. Therefore,
only the branch CNN remains fixed to the $50\times50$ distance-stack input,
whereas the remaining operator components can be queried on the finer
$100\times100$ grid.

In the three-point-source case with the fixed smooth velocity
$F(x_1,x_2)=2x_1^2+1$, training and evaluating the model on the
$50\times50$ grid gives a relative $L^2$ error of $4.12\cdot 10^{-2}$. When the same
trained model is evaluated on the original $100\times100$ grid, while keeping
the branch input at $50\times50$, the relative $L^2$ error is $4.27\cdot 10^{-2}$.
Thus, the error changes only marginally despite doubling the evaluation
resolution in each spatial direction. This indicates that the learned
operator is not tightly tied to the training discretization at the
output-evaluation stage: the branch CNN provides a coarse but stable encoding
of the source geometry, while the following FNO blocks and lower-envelope output transform reconstruct the travel-time field on the finer query grid.

\section{Conclusion}\label{sec:conclusion}
We proposed Lower-envelope Eikonal Neural Operator (LENO), a physics-informed neural operator for the eikonal equations with general source sets and heterogeneous velocity fields. The key idea is to encode the variational minimum principle of the viscosity solution directly into the architecture: the model constructs source-wise candidate fields and takes their pointwise
minimum. This lower-envelope transform enforces the source condition $u|_\Gamma=0$ exactly and represents the multi-source structure of the solution without a soft source penalty.

The model is trained without precomputed solution data, using a Godunov upwind residual as the physics-informed training signal. Numerical experiments show that the same architecture handles multiple point sources, arbitrary open or closed curve sources, sample-dependent piecewise velocity fields, and OpenFWI CurveVel-A and Style-A velocity models with surface sources. Ablation results confirm that the velocity-augmented trunk features and the lower-envelope output transform are essential in realistic heterogeneous media, while the resolution-transfer experiment shows that the trained operator can be queried on a finer output grid with only marginal degradation. In addition to its accuracy, LENO exhibits useful operator-level properties. Its single-instance inference time is comparable to that of classical Fast Marching and Sweeping Methods, while batched evaluation substantially reduces the total runtime when many source-velocity configurations must be queried. 

Future work will extend the framework to three-dimensional problems,
unstructured discretizations, and inverse settings in which the source set or velocity field is inferred from partial travel-time observations. Moreover, because LENO represents the solution map by a differentiable neural operator, $\mathcal{S}_\theta:(\Gamma,F)\mapsto u$, it is naturally compatible with gradient-based downstream tasks. This differentiability is particularly relevant for inverse problems such as traveltime tomography and source localization, where classical solver-based approaches often require separately
derived adjoint-state formulations~\cite{tong2021adjoint}.

\section*{Declaration of generative AI and AI-assisted technologies in the manuscript preparation process}
During the preparation of this work the authors used CLAUDE in order to improve English. After using this tool/service, the authors reviewed and edited the content as needed and take full responsibility for the content of the published article.
\printcredits

\section*{Acknowledgements}
Hwijae Son was supported by Basic Science Research Program through the National Research Foundation of Korea (NRF) funded by the Ministry of Education (RS-2024-00462755, RS-2026-25491905). Sung Woong Cho was supported by the National Research Foundation of Korea (NRF) grant funded by the Ministry of Science and ICT (KR) (RS-2024-00462755) and was supported by Inha University.

\section*{Code Availability}
The original code that supports the findings of this study is available from the corresponding author upon reasonable request.

\bibliographystyle{cas-model2-names}
\bibliography{references}







\end{document}